\documentclass[12pt, a4paper]{amsart}

\usepackage{amsmath, amssymb}

\usepackage[a4paper]{geometry}
\usepackage{enumerate}

\usepackage{graphicx}
\usepackage{wrapfig}

\newtheorem{theorem}{Theorem}

\theoremstyle{remark}

\newcommand{\C}{\ensuremath{\mathbb{C}}}
\newcommand{\R}{\ensuremath{\mathbb{R}}}

\newcommand{\cal}[1]{\ensuremath{\mathcal{#1}}}

\DeclareMathOperator{\spann}{span}

\begin{document}
	
\title[Inhomogeneous ruled hypersurface]{Existence and uniqueness of inhomogeneous ruled hypersurfaces with shape operator\\ of constant norm in the complex hyperbolic space}

\author[M.~Dom\'{\i}nguez-V\'{a}zquez]{Miguel Dom\'{\i}nguez-V\'{a}zquez}
\address{Department of Mathematics, University of Santiago de Compostela, Spain.}
\email{miguel.dominguez@usc.es}
\author[O.~P\'erez-Barral]{Olga P\'erez-Barral}
\address{Department of Mathematics, University of Santiago de Compostela, Spain.}
\email{olgaperez.barral@usc.es}

\thanks{The authors have been supported by projects MTM2016-75897-P, PID2019-105138GB-C21 (AEI/FEDER, Spain) and ED431C 2019/10, ED431F 2020/04 (Xunta de Galicia, Spain). The first author acknowledges support of the	Ram\'{o}n y Cajal program of the Spanish State Research Agency.}

\subjclass[2010]{53B25, 53C42, 53C55}


\begin{abstract}
	We complete the classification of ruled real hypersurfaces with shape operator of constant norm in nonflat complex space forms by showing the existence of a unique inhomogeneous example in the complex hyperbolic space. 
\end{abstract}

\keywords{Complex hyperbolic space, complex space form, ruled hypersurface, constant norm, shape operator}

\maketitle

\section{Introduction}

A real hypersurface of a complex space form is called ruled if it is locally foliated by totally geodesic complex hypersurfaces of the ambient space. The aim of this note is to prove the following classification result. In what follows,  $\C H^n$ denotes the $n$-dimensional complex hyperbolic space of constant holomorphic sectional curvature $c<0$.

\begin{theorem}\label{theorem:ruled:II:new:example}
	Let $M$ be a ruled real hypersurface in a nonflat complex space form. The shape operator of $M$ has constant norm if and only if $M$ is an open part of:
	\begin{enumerate}[{\rm (i)}]
		\item A Lohnherr hypersurface of $\mathbb{C}H^{n}$, or
		\item The ruled real hypersurface of $\C H^n$ constructed by attaching totally geodesic $\mathbb{C}H^{n-1}$ perpendicularly to an entire circle of curvature $\sqrt{-c/2}$ in a totally geodesic $\mathbb{C}H^{1}$.
	\end{enumerate}
\end{theorem}

The Lohnherr hypersurface mentioned in  Theorem~\ref{theorem:ruled:II:new:example}~(i) is a homogeneous hypersurface and, indeed, it is the only ruled real hypersurface with constant principal curvatures of a nonflat complex space form~\cite{LR}.  It has been investigated from different viewpoints in the literature; see~\cite{DVPB:jgp} for references. Recently, the second author proved the inexistence of ruled hypersurfaces with shape operator of constant norm in complex projective spaces $\C P^n$, and showed that any possible example in $\C H^n$ should be strongly $2$-Hopf~\cite{PB:results}. 

The novelty of this article is to show that there is exactly one inhomogeneous complete ruled hypersurface with shape operator of constant norm in $\C H^n$, and it is constructed as in Theorem~\ref{theorem:ruled:II:new:example}~(ii). Thus, Theorem~\ref{theorem:ruled:II:new:example} constitutes one more contribution to the general problem of identifying and characterizing the ``simplest" ruled hypersurfaces, cf.~\cite[p.~446]{CR}. In this sense, we remark that
two recent papers investigating two other geometric properties, namely the constancy of the mean curvature~\cite{DVPB:jgp} or of the scalar curvature~\cite{KMT:gd}, either characterize well-known examples (the minimal ruled hypersurfaces~\cite{LR}) in the former case, or a rather broad collection of examples (parametrized by lightlike curves in an indefinite complex projective space) in the latter. Finally, we notice that the  the constancy of the norm of the shape operator is a classical property in Differential Geometry, insofar as it arises naturally in important problems, such as the longstanding Chern's conjecture~\cite{AB:duke}.

\section{Proof of Theorem~\ref{theorem:ruled:II:new:example}}\label{section:ruled:II}

Let $M$ be a ruled real hypersurface of a nonflat complex space form whose shape operator $\cal{S}$ has constant norm. By \cite[Theorem~1]{PB:results} we can assume that the ambient space is a complex hyperbolic space $\C H^{n}$, $n\geq 2$, since such a hypersurface $M$ cannot exist in complex projective spaces. We denote by $J$ the complex structure of $\C H^n$, and let $\xi$ be a (locally defined) unit normal vector field to $M$. By definition, the distribution $(J\xi)^\perp$ of tangent vectors to $M$ which are orthogonal to the Hopf field $J\xi$ is autoparallel, and its integral submanifolds are open subsets of totally geodesic $\C H^{n-1}$. The shape operator of $M$ satisfies $\cal{S}J\xi=\nu J\xi+\mu JA$, $\cal{S}JA=\mu J\xi$ and $\cal{S}X=0$ for any $X\in TM$ perpendicular to $\spann\{J\xi,JA\}$, where $A$ is a (locally defined) unit vector field in the $0$-principal curvature distribution of $M$. See~\cite[\S2]{ABM}, \cite[\S8.5.1]{CR} or \cite[Propositions~3.1 and 3.2]{DVPB:jgp} for details.

Now we recall from~\cite[Theorem~1]{PB:results} that $M$ is a strongly $2$-Hopf real hypersurface. This means that the smallest $\cal{S}$-invariant distribution containing the Hopf vector field, $\cal{D}=\spann\{J\xi,JA\}$, has rank 2, is integrable and the spectrum of $\cal{S}\vert_\cal{D}$ is constant along the integral submanifolds of $\cal{D}$. Moreover, by the discussion in~\cite[end of \S3.1]{PB:results} and the fact that any ruled hypersurface with constant principal curvatures is an open part of a Lohnherr hypersurface~\cite[Remark~5]{LR}, whose nonzero principal curvatures are $\pm\sqrt{-c}/2$, we have that $M$ satisfies $|\cal{S}|^2=-c/2$. Since $|\cal{S}|^2=2\mu^2+\nu^2$, we obtain $|\mu|\leq \sqrt{-c}/2$ on $M$.

Assume that there is a point $p\in M$ with $|\mu(p)|=\sqrt{-c}/2$. On the one hand, since $M$ is strongly $2$-Hopf, the eigenvalues of $\cal{S}\vert_\cal{D}$, or equivalently the functions $\mu$ and $\nu$, are constant along the integral submanifolds of $\cal{D}=\spann\{J\xi,JA\}$. On the other hand, by~\cite[Corollary~3]{LR} $\mu=\pm\sqrt{-c}/2$ is constant along the integral submanifold of $(J\xi)^\perp$ containing $p$. All in all, we deduce that $\mu=\pm\sqrt{-c}/2$ is constant on $M$, and since $2\mu^2+\nu^2=-c/2$, then $\nu=0$ on $M$. Thus, $M$ has constant principal curvatures, and hence, $M$ is an open part of a Lohnherr hypersurface, which corresponds to case (i) in Theorem~\ref{theorem:ruled:II:new:example}.

From now on we will assume that $|\mu|<\sqrt{-c}/2$ on $M$. Then, the extension $\widetilde{M}=\bigcup_{p\in M}\exp_p(J\xi_p)^\perp$ of $M$, given by the union of every totally geodesic $\C H^{n-1}$ that contains an integral submanifold of $(J\xi)^\perp$, is a ruled hypersurface without singular points by~\cite[Theorem~2]{LR}. We will still denote by $\xi$ and $J\xi$ the unit normal and Hopf vector fields of $\widetilde{M}$, respectively. By the standard theory of ruled hypersurfaces (see~\cite[Theorem~3 along with (26) and (29)]{LR}), the functions $\nu$ and $\mu^2$ are real analytic when restricted to any geodesic contained in a ruling of $\widetilde{M}$. Since $2\mu^2+\nu^2=-c/2$ on the open subset $M$ of $\widetilde{M}$, the same relation holds on the whole $\widetilde{M}$. Then, $\widetilde{M}$ is also strongly $2$-Hopf. Moreover, since $\widetilde{M}$ is smooth and with complete rulings, it follows from~\cite[Corollary~3]{LR} (cf.~\cite[Lemma~1]{ABM}) 
that there is a point $p\in\widetilde{M}$ where $\mu(p)=0$ and, thus, $\nu(p)=\pm\sqrt{-c/2}$. But as $\widetilde{M}$ is strongly $2$-Hopf, both $\nu$ and $\mu$ are constant along the integral curves of the Hopf vector field  of $\widetilde{M}$. Thus, if $\gamma$ is the integral curve of $J\xi$ with $\gamma(0)=p$, then $\nu(\gamma(t))=\pm\sqrt{-c/2}$ and $\mu(\gamma(t))=0$, for all $t\in\R$ where $\gamma$ is defined.
Moreover, a basic identity for ruled hypersurfaces ensures that 
\begin{equation}\label{eq:nablaJxi}
\bar{\nabla}_{J\xi}J\xi=\mu A+\nu\xi,
\end{equation} 
where $\bar{\nabla}$ is the Levi-Civita connection of the ambient space $\C H^n$. Then, the integral curve $\gamma$ of $J\xi$ with $\gamma(0)=p$ satisfies the differential equations

\begin{align}\label{eq:curve}
	&\bar{\nabla}_{\dot{\gamma}}\dot{\gamma}=\mp\sqrt{\frac{-c}{2}} J\dot{\gamma}, &&\bar{\nabla}_{\dot{\gamma}}\bar{\nabla}_{\dot{\gamma}}\dot{\gamma}=\frac{c}{2 } \dot{\gamma},
\end{align} 
which means that $\gamma$  is a circle of curvature $\sqrt{-c/2}$ inside a totally geodesic complex hyperbolic line $\mathbb{C}H^{1}$ (see~\cite[Section 2]{ABM}). Thus, we have shown that $M$ is an open part of the ruled hypersurface $\widetilde{M}$ constructed by attaching totally geodesic $\C H^{n-1}$ perpendicularly to a circle $\gamma$ of curvature $\sqrt{-c/2}$ inside a totally geodesic $\mathbb{C}H^{1}$, which corresponds to case (ii) in Theorem~\ref{theorem:ruled:II:new:example}.

In order to finish the proof of Theorem~\ref{theorem:ruled:II:new:example} we have to justify that the ruled hypersurface constructed as in (ii) is smooth and its shape operator has constant norm. Let then $\gamma\colon \R\to\C H^n$  be a complete circle of curvature $\sqrt{-c/2}$ in  a totally geodesic  $\mathbb{C}H^{1}$; note that, as $\sqrt{-c/2}<\sqrt{-c}$, $\gamma$ is an equidistant curve to a geodesic in the totally geodesic $\C H^1$. Let $M=\bigcup_{t\in\R}\exp_{\gamma(t)}(\C\dot{\gamma})^\perp$ be the corresponding ruled hypersurface.  Since $M$ and $\gamma$ satisfy~\eqref{eq:nablaJxi} and~\eqref{eq:curve}, respectively, and  $\dot{\gamma}=J\xi$ along $\gamma$ (maybe changing the sign of $\xi$), we get $\mu\circ\gamma=0$ and $\nu\circ\gamma=\sqrt{-c/2}$. Appealing again to~\cite[Theorems~2 and~3]{LR} we deduce that $M$ is an immersed hypersurface of $\C H^n$ and the functions $\mu$ and $\nu$ are given by
\[
\mu(\sigma(r))=\pm\frac{\sqrt{-c}}{2}\tanh\left(\frac{\sqrt{-c}}{2}r\right),\qquad\qquad \nu(\sigma(r))=\sqrt{\frac{-c}{2}}\,\mathrm{sech}\left(\frac{\sqrt{-c}}{2}r\right),
\]
where $\sigma(r)=\exp_{\gamma(t)}(rX)$ is unit speed geodesic in the ruling of $M$ through $\gamma(t)$, for any unit $X\in(\C\dot{\gamma}(t))^\perp$ and any $t\in\R$. Thus, it follows that the squared norm of the shape operator of $M$ satisfies $|\cal{S}|^2=2\mu^2+\nu^2=-c/2$, as we wanted to show. Finally, note that $M$ is closed and embedded in $\C H^n$, as follows easily from the fact that the totally geodesic complex hyperplanes $\C H^{n-1}$ that are perpendicular to a totally geodesic $\C H^1$ define a smooth foliation of $\C H^n$.

\end{document}